\documentclass[12pt]{article}  

\usepackage{amsmath}
\usepackage{amssymb}
\usepackage{cite}
\usepackage[usenames]{color}
\usepackage{graphicx}          
\usepackage[dvips]{epsfig}    
\newtheorem{definition}{Definition}
\newtheorem{theorem}{Theorem}
\newtheorem{corollary}{Corollary}
\newtheorem{remark}{Remark}
\newtheorem{problem}{Problem}
\newtheorem{lemma}{Lemma}

\title{\LARGE \bf Stabilization of  non-admissible curves for a class of nonholonomic systems
\thanks{
This work was  supported in part by the German Research Foundation (GR 5293/1-1), NAS of Ukraine (budget program KPKBK 6541230), and the State Fund for Fundamental Research of Ukraine (F75/27190)\newline
$^{1}$Institute of Mathematics, University of W\"{u}rzburg,       97074 W\"{u}rzburg, Germany
        {\tt\small viktoriia.grushkovska@mathematik.uni-wuerzburg.de}
        \newline
$^{2}$Max Planck Institute for Dynamics of Complex Technical Systems, 39106 Magdeburg, Germany
{\tt\small zuyev@mpi-magdeburg.mpg.de}\newline
$^{3}$Institute of Applied Mathematics and Mechanics, National Academy of Sciences of Ukraine, 841116 Sloviansk, Ukraine
}
}

\author{Victoria Grushkovskaya$^{1,3}$ and Alexander Zuyev$^{2,3}$}
\date{}

\begin{document}

\maketitle
\thispagestyle{empty}

\begin{abstract}
The problem of tracking an arbitrary curve in the state space is considered for underactuated driftless control-affine systems.
This problem is formulated as the stabilization of a time-varying  family of sets associated with a neighborhood of the reference curve.
An explicit control design scheme is proposed for the class of controllable systems whose degree of nonholonomy is equal to 1.
It is shown that the trajectories of the closed-loop system converge exponentially to any given neighborhood of the reference curve provided that the solutions are defined in the sense of sampling.
This convergence property is also illustrated numerically by several examples of nonholonomic systems of degrees 1 and 2.
\end{abstract}

\section{INTRODUCTION}
In this paper, we consider a class of driftless control systems of the form
\begin{equation}
\dot x = \sum_{i=1}^m u_i f_i (x),\quad x\in \mathbb R^n,\;  u\in \mathbb R^m,\;m<n, \; f_i \in C^2(\mathbb R^n),
\label{Sigma}
\end{equation}
where $x=(x_1,\dots,x_n)^\top$ is the state and $u=(u_1,\dots,u_m)^\top$ is the control.
The stabilization of such systems has been the subject of numerous studies over the last few decades, and  many important results have been obtained in this area.
In particular, it follows from the famous result of R.W.~Brockett~\cite{Bro81} that the trivial equilibrium of~\eqref{Sigma} is not stabilizable by a regular time-invariant feedback law if the vectors $f_1(0)$, $f_2(0)$, ..., $f_m(0)$ are linearly independent.
Despite the significant progress in the development of control algorithms to stabilize the solution $x=0$ of system~\eqref{Sigma} (see, e.g.,~\cite{Ast94,Bloch92,Cor92,Kol95,Pan11,Pom92,Sar17,Tian02,Zu16,ZG16}, and references therein), the  stabilization of nonholonomic systems to a given curve  remains a challenging problem.
This issue can be formulated as the trajectory tracking problem. In many papers, this problem has been addressed under the assumption that the trajectory is admissible, i.e. satisfies the system equations with some control inputs~\cite{Ai05,Ali16,Dandr95,Dong99,Fl95,Mag17,Walsh92,Wang15,Yu15}.
 Since the number of controls $m$ may be significantly smaller than the dimension of the state space $n$, not every path in the state space is admissible for system~\eqref{Sigma}.
However, in many applied problems, it is important to stabilize system~\eqref{Sigma} along an \emph{arbitrary} curve, which is not necessarily admissible. As it is mentioned in~\cite{MS08b}, although it is not possible to asymptotically stabilize nonholonomic systems to non-admissible curves because of the non-vanishing tracking error, the practical stabilization can be achieved. It has to be noted that such problem has been addressed only for particular classes of systems, e.g., for unicycle and car-like systems~\cite{MS08b,GMZME18,Rav18}


This paper deals with rather general formulation of the stabilization problem with non-admissible reference curves.
The main contribution of our paper is twofold. First, we introduce a class of control functions for the first degree nonholonomic systems, which allows  stabilizing the system in a prescribed neighborhood of an arbitrary (not necessarily admissible) curve. We also show how the obtained results can be extended to higher degree nonholonomic systems. The proposed feedback design scheme is based on the approach introduced in~\cite{Zu16,ZG17,GZ18} for the stabilization and motion planning of nonholonomic systems. However, it has to be noted that the results of these papers cannot be directly applied for the stabilization of non-admissible curves. Second, we characterize stability properties of system~\eqref{Sigma} with the proposed controls in terms of families of sets. Note that the concept of stability of families of sets was used previously in~\cite{La02,GDEZ17,GMZME18} for non-autonomous system admitting a Lyapunov function.  In the present paper, we do not assume the existence of a control Lyapunov function and define solutions of the closed-loop system in the sense of sampling.

The rest of the paper is organized as follows. In the remainder of this section, we introduce some basic notations, recall the notion of stability of sets, and give a precise problem statement.
The main result will be proved in Section~II and  illustrated with some examples in Section~III.
\subsection{Notations and definitions}
To generate attractive control strategies for system~\eqref{Sigma} in a neighborhood of a given curve $\Gamma=\{\gamma(t)\}_{t\ge 0}\subset \mathbb R^n$, we will follow the idea of~\cite{Zu16} and define solutions of the corresponding closed-loop system in the sense of sampling.
With a slight abuse of notation, we will also identify the curve $\Gamma=\{\gamma(t)\}_{t\ge 0}$ with the map $\gamma:\mathbb R^+\to \mathbb R^n$,  $\mathbb R^+= [0,+\infty)$.
For a given $\varepsilon>0$, we consider the partition $\pi_\varepsilon$ of $\mathbb R^+$ into intervals
$
I_j=[t_j,t_{j+1}),\;t_j=\varepsilon j, \quad j=0,1,2,\dots \; .
$
\begin{definition}
Assume given a curve $\gamma:\mathbb R^+ \to \mathbb R^n$, a feedback law $h: \mathbb R^+ \times \mathbb R^n \times \mathbb R^n  \to \mathbb R^m$, and an~$\varepsilon>0$.
A $\pi_\varepsilon$-solution of~\eqref{Sigma} corresponding to $x^0\in \mathbb R^n$ and $u=h(t,x,\gamma)$ is an absolutely continuous function  $x(t)\in \mathbb R^n$, defined for $t\in[0,+\infty)$, such that  $x(0)=x^0$ and, for each $j=0, 1, 2, \dots$,
$$
\dot x(t)=\sum_{i=1}^mh_i(t,x(t_j),\gamma(t_j))\big)f_i(x(t)), \quad t\in I_j=[t_j,t_{j+1}).
$$

\end{definition}

  For  $f,g:\mathbb R^n\to\mathbb R^n $, $x\in\mathbb R^n$, we denote the Lie derivative as
 {$ L_gf(x)=\lim\limits_{s\to0}\tfrac{f(x+sg(x))-f(x)}{s}$}, and  $[f,g](x)= L_fg(x)- L_gf(x)$ is the Lie bracket.
 Throughout this paper, $\|a\|$ stands for the Euclidean norm of a vector $a\in\mathbb R^n$, and the norm of an $n\times n$-matrix $\cal F$ is defined as $\|{\cal F}\|=\sup_{\|y\|=1}\|{\cal F}y\|$.

\subsection{Stability of a family of sets}
To characterize the asymptotic behavior of trajectories of system~\eqref{Sigma},
we will extend the concept of stability of a family of sets to the case of $\pi_\varepsilon$-solutions.
This concept has been developed, e.g., in~\cite{La02} for non-autonomous differential equations and applied to control problems under the classical definition of solutions in~\cite{GDEZ17,GMZME18}.
   Let
$\{\mathcal S_t\}_{t\ge 0}$ be a one-parameter family of non-empty subsets of $\mathbb R^n$.
For a $\delta>0$, we denote the $\delta$-neighborhood of the set
$\mathcal S_{t}$ at time $t$ as
 $
 B_{\delta}(\mathcal S_{t}){=}\bigcup_{y{\in}\mathcal S_{t}}\{x{\in}\mathbb R^n:\|x{-}y\|{<}\delta\},
 $
  The distance from a point $x\in \mathbb R^n$ to a set $\mathcal S_{t}\subset \mathbb R^n$ is denoted as ${\rm dist}(x,\mathcal S_{t})=\inf_{y\in \mathcal S_{t}}\|x-y\|$.
  Assume given a curve $\gamma:\mathbb R^+ \to \mathbb R^n$, a time-varying feedback law $h: \mathbb R^+ \times \mathbb R^n \times \mathbb R^n  \to \mathbb R^m$, and a sampling parameter~$\varepsilon>0$. The basic stability definition that we exploit in this paper is as follows.
\begin{definition}
A one-parametric family of  sets   $\{\mathcal S_{t}\}_{t\ge0}$ is said to be \emph{exponentially stable} for the closed-loop system~\eqref{Sigma} with $u=h(t,x,\gamma)$ \emph{in the sense of} $\pi_\varepsilon$-\emph{solutions} if
   there exist $\hat \delta,\lambda>0$ such that, for any $x^0{\in} B_{\hat\delta}(\mathcal S_{0})$, the corresponding
   $\pi_\varepsilon$-solution of~\eqref{Sigma} satisfies ${\rm dist}(x(t),\mathcal S_t)\le C e^{-\lambda t}$ for all $t \ge 0$ with some $C=C(x^0)$.%
  If the above exponential decay property holds for every $\hat\delta{>}0$, then the family of sets $\{\mathcal S_{t}\}_{t\ge0}$ is called  \emph{globally exponentially stable}
  in the sense of $\pi_\varepsilon$-solutions.
 \end{definition}
\subsection{Problem statement}~\label{sec_problem}
Using the notion of stability of a family of sets, it is convenient to formulate the control design problem under consideration as follows:
\begin{problem}
  Given a curve $\gamma:\mathbb R^+\to\mathbb R^n$ and a constant $\rho>0$, the goal is to find a time-varying feedback law $h: \mathbb R^+ \times \mathbb R^n \times \mathbb R^n\to \mathbb R^m$ such that the family of sets
\begin{align}\label{set}
    \{\mathcal S_t^\rho\}_{t\ge0}=\left\{\mathcal S_t^\rho=B_\rho(\gamma(t))\}\right\}_{t\ge0}
\end{align}
is exponentially stable for the closed-loop system~\eqref{Sigma} with $u=h(t,x,\gamma)$ in the sense of Definition~2.
\end{problem}

We will propose a solution to the above problem with a $C^1$-curve $\gamma:\mathbb R^+\to\mathbb R^n$ for the nonholonomic systems of degree one, i.e.,
we assume that there is an $r>\rho$ such that the following rank condition holds in $D=\bigcup_{t\ge0}B_r(\gamma(t))$:
\begin{equation}\label{rank}
 {\rm span}\big\{f_{i}(x), [f_{j_1},f_{j_2}](x):\,i{\in}S_1,(j_1,j_2){\in} S_2\big\}=\mathbb{R}^n
\end{equation}
for all $x\in D$, with some sets of indices $S_1\subseteq \{1,2,...,m\}$,  $S_2\subseteq \{1,2,...,m\}^2$ such that $|S_1|+|S_2|=n$.

\section{MAIN RESULTS}
\subsection{Control design}
To solve Problem~1, we extend the control design approach proposed in~\cite{GZ18}. Namely, we use a family of trigonometric polynomials with state-dependent coefficients chosen in such a way that the trajectory of system~\eqref{Sigma} approximate the gradient flow of a  time-invariant Lyapunov function.
In this paper, the corresponding Lyapunov function is time-varying, so we allow the above mentioned coefficients  to depend on time.
We define the control functions in the following way:
\begin{align}
  &u_i^\varepsilon(t,x,\gamma)=\sum_{j\in S_1}\delta_{i j} a_{j}(x,\gamma)\nonumber\\
  &+\sqrt{\frac{4\pi}{\varepsilon}}\sum_{(j_1,j_2)\in S_2}{\sqrt{\kappa_{j_1j_2}|a_{j_1,j_2}(x,\gamma)|}}\Big(\delta_{ij_1}\cos{\frac{2\pi \kappa_{j_1j_2}}{\varepsilon}}t\nonumber\\
  &\quad{+}\delta_{ij_2}{\rm sign}(a_{j_1,j_2}(x,\gamma))\sin{\frac{2\pi \kappa_{j_1j_2}}{\varepsilon}}t\Big),\; i=1,2,...,m.\label{cont}
\end{align}
Here $\delta_{ij}$ is the Kronecker delta, $\kappa_{j_1j_2}{\in}\mathbb N$ are pairwise distinct, and
$$\Big((a_{j}(x,\gamma))_{j \in S_1}\  ( a_{j_1j_2}(x,\gamma))_{(j_1,j_2)\in S_2}\Big)^\top=a(x,\gamma),$$
where
\begin{equation}\label{a}
a(x,\gamma)=- \alpha \mathcal F^{-1}(x) (x-\gamma),
\end{equation}
with
$
  \mathcal F(x){= }\Big(\big(f_{j}(x)\big)_{j\in S_1}\ \ \big([f_{j_1},f_{j_2}](x)\big)_{(j_1,j_2)\in S_2}\Big)
$
and  $\alpha{>}0$.
Note that~\eqref{rank} implies that $\mathcal F(x)$ is nonsingular in $ D$.
\subsection{Stability analysis}
The main result of this paper is as follows.
\begin{theorem}\label{main}
\emph{
Let $\gamma\in C^1(\mathbb R^+;\mathbb R^n)$, 
$r> 0$, $\mu>0$, and $\nu\ge 0$ be such that the matrix $\mathcal F(x)$  is nonsingular in
$\displaystyle D=\bigcup_{t\ge0}B_r(\gamma(t))$,
$
f_i(x)$, $L_{f_j} f_i(x)$, $L_{f_k} L_{f_j} f_i(x)$ are bounded  in~$D$ ($i,j,k=\overline{1,m}$), $\|\mathcal F^{-1}(x)\|\le \mu$ for all $x\in D$,
and $\|\dot \gamma(t)\|\le\nu$ for all $t\ge 0$.
Then, for any $\rho\in(0,r)$, there exists an $\hat \varepsilon>0$ such that the family of sets~\eqref{set}
is exponentially stable for system~\eqref{Sigma} with the controls $u_i=u_i^\varepsilon$ defined by~\eqref{cont}--\eqref{a} with any $\varepsilon\in (0,\hat \varepsilon)$ and $\alpha>\frac{\nu}{\rho}$ in the sense of Definition~2.}
\end{theorem}
The proof of this theorem is given in the Appendix.

The next corollary follows from the proof of Theorem~1.
\begin{corollary}
\emph{Let the conditions of Theorem~\ref{main} be satisfied, and let $\|\dot\gamma(t)\|\to 0$ as $t\to+\infty$. Then there is a $\hat \delta>0$ such that $\|x(t)-\gamma(t)\|\to 0$ as $t\to +\infty$, provided that $\|x(0)-\gamma(0)\|<\hat \delta$ and the solutions of the closed-loop system~\eqref{Sigma},~\eqref{cont}--\eqref{a} are defined in the sense of Definition~1.
}\end{corollary}
%
%
  %

Let us emphasize that, in contrast to many other results on stability of non-autonomous systems (e.g.,~\cite{Khalil}), we do not require the boundedness of $\gamma(t)$ in general.


\section{EXAMPLES}
In this section, we consider some examples illustrating Theorem~\ref{main} and discuss the possibility of extending the above results to systems with a higher degree of nonholonomy.
\begin{figure*}[h!]
 \begin{minipage}{0.33\linewidth}
  \includegraphics[width=1\linewidth]{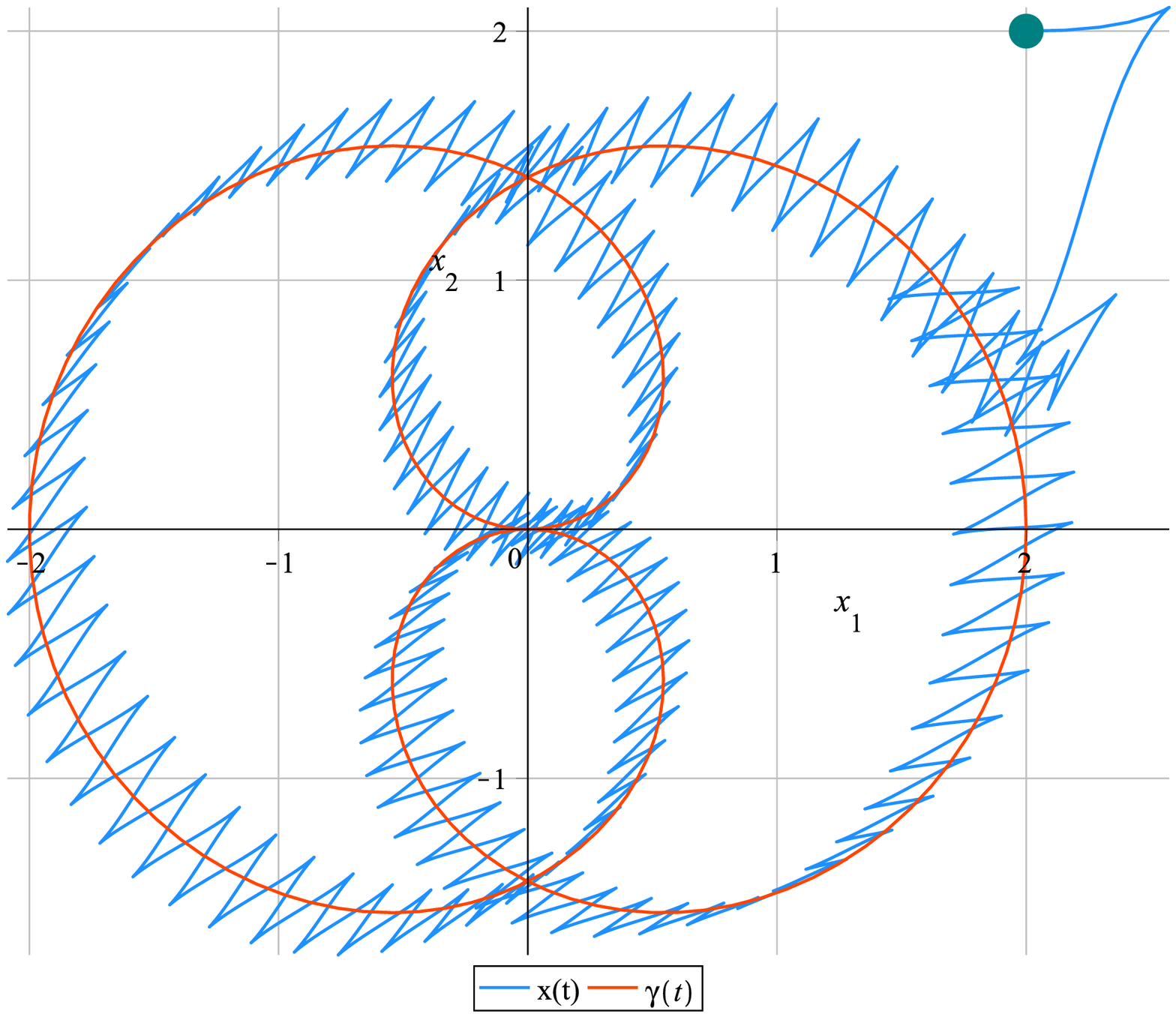}
\includegraphics[width=1\linewidth]{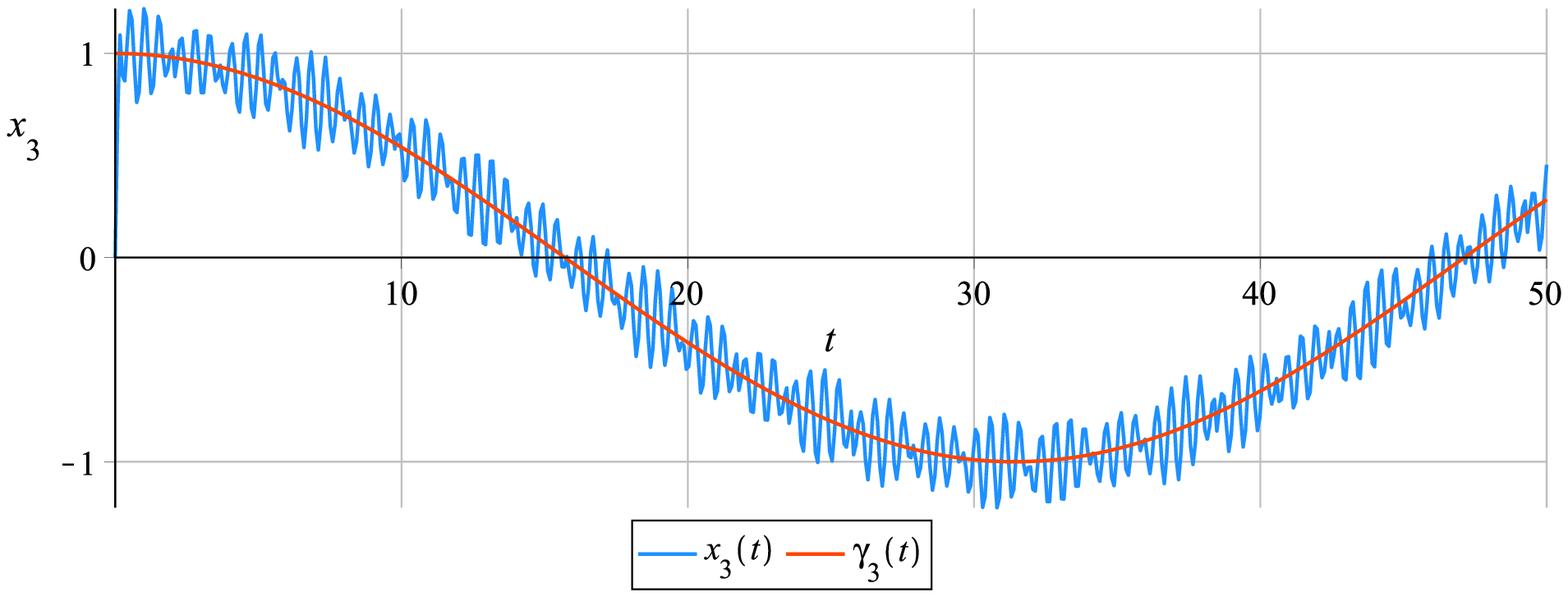}
 \end{minipage}\hfill
  \begin{minipage}{0.33\linewidth}
  \includegraphics[width=1\linewidth]{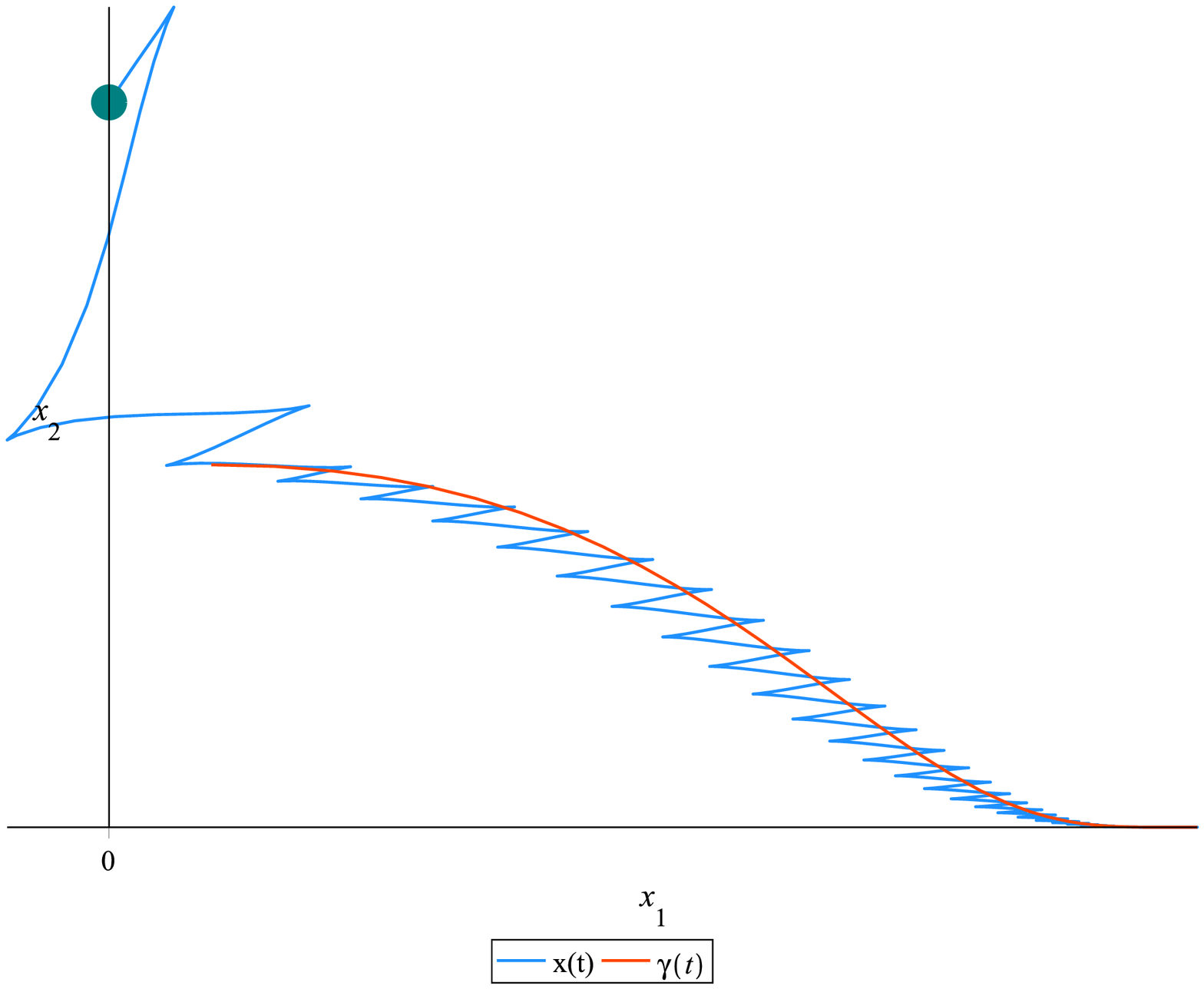}
\includegraphics[width=1\linewidth]{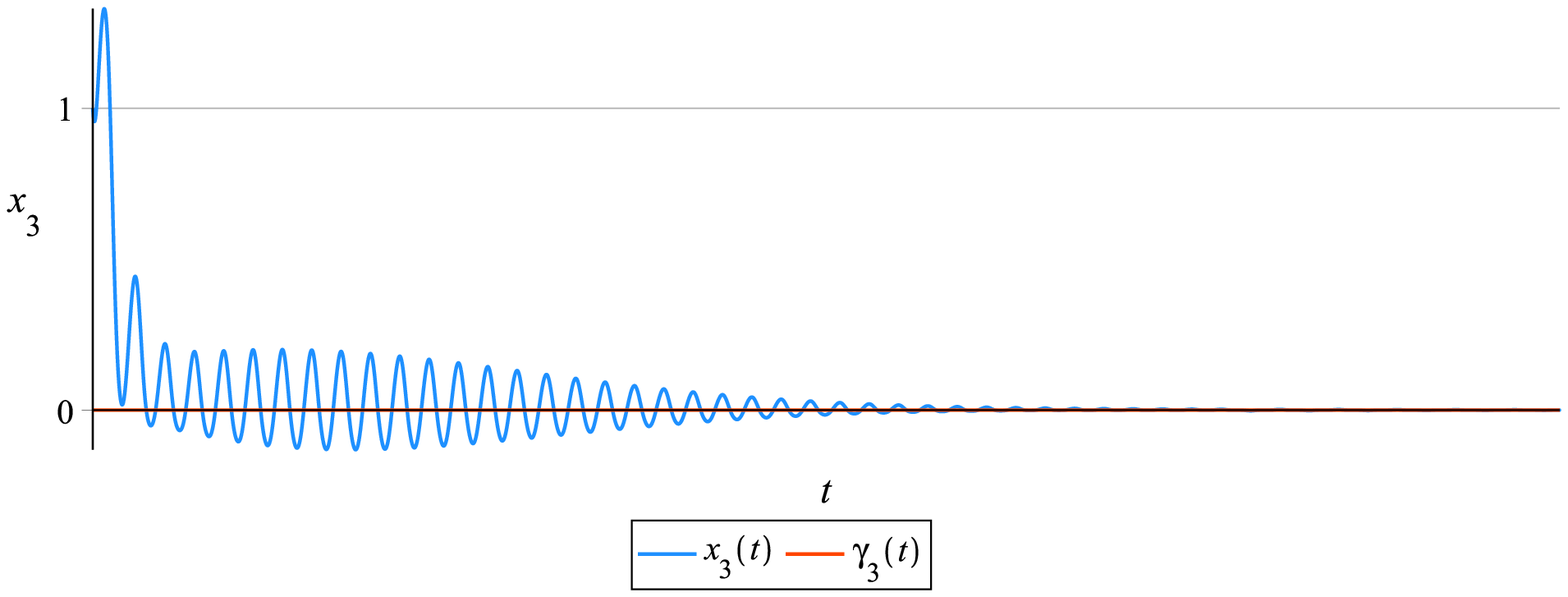}
 \end{minipage}\hfill
 \begin{minipage}{0.33\linewidth}
  \includegraphics[width=0.99\linewidth]{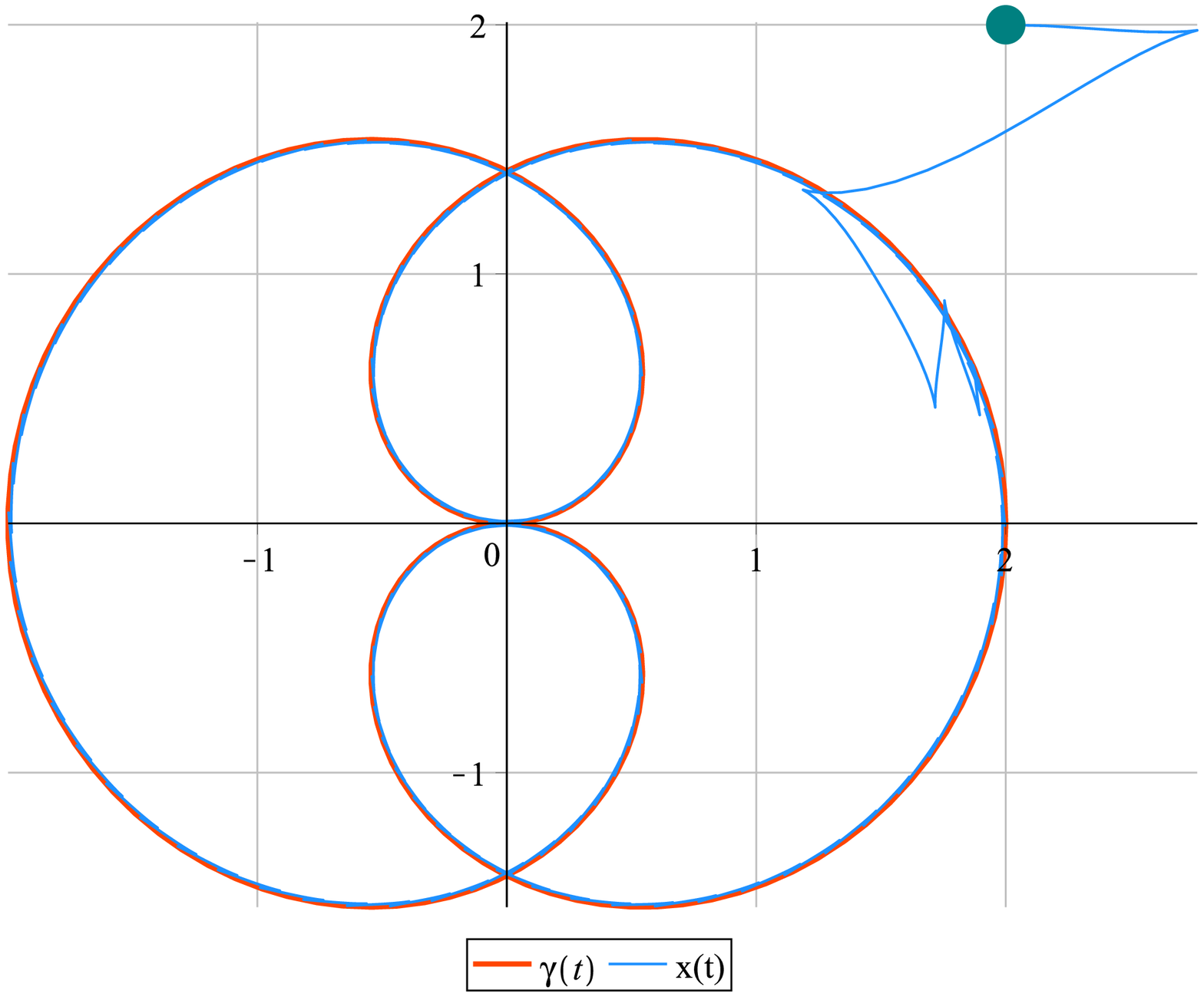}
\includegraphics[width=1\linewidth]{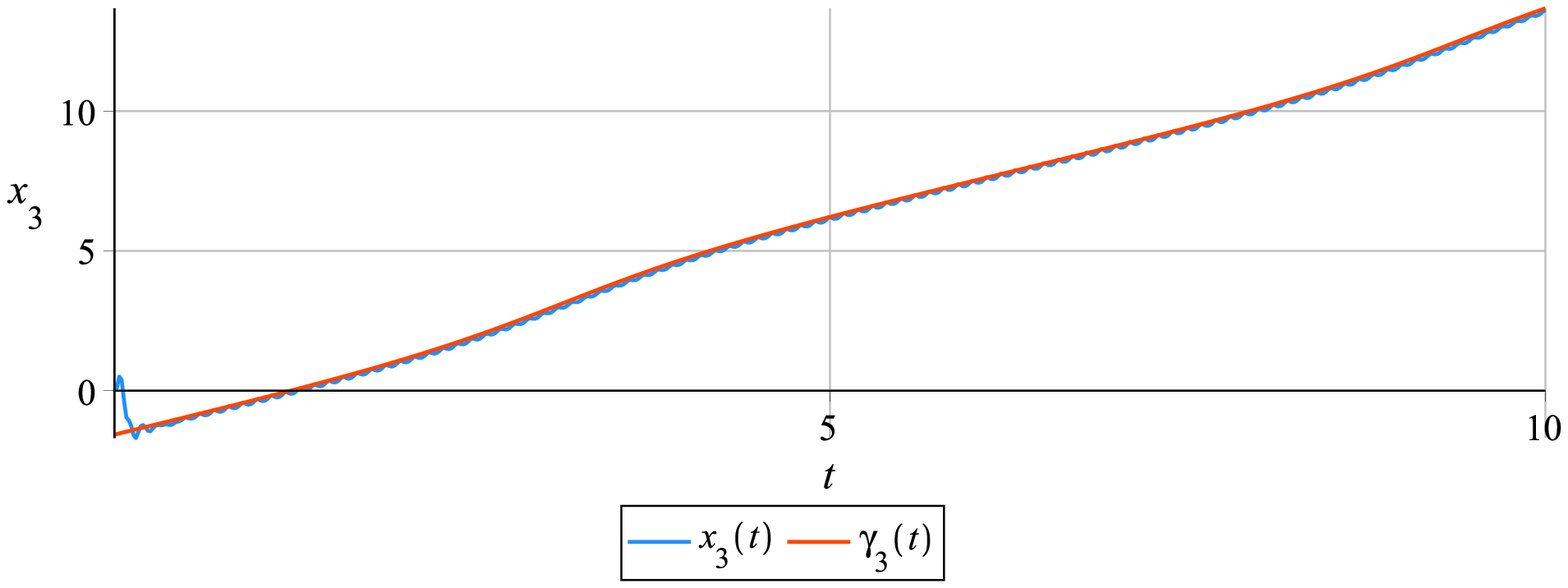}\label{uni_feas}
 \end{minipage}
  \caption{Trajectories of system~\eqref{ex_uni} with controls~\eqref{cont_uni} ($\alpha=15$, $\varepsilon=0.1$) and the curves $\gamma^{(1)}$ (left),  $\gamma^{(2)}$ (middle),  $\gamma^{(3)}$ (right). }
\end{figure*}
\subsection{Unicycle}
As the first example, consider the equations of motion of the unicycle:
\vskip-2ex
\begin{equation}\label{ex_uni}
\dot x_1=u_1\cos x_3,\ \dot x_2=u_1\sin x_3,\ \dot x_3=u_2,
\end{equation}
where $(x_1,x_2)$ are the coordinates of the contact point of the unicycle, $x_3$ is the angle between the wheel and the $x_1$-axis, $u_1$ and $u_2$ control the forward and the angular velocity, respectively. Denote $f_1(x)=\big(\cos (x_3),\sin (x_3), 0\big)^\top$, $f_2(x)=\big(0,0,1\big)^\top$. Then the rank condition~\eqref{rank} is satisfied for all $x\in \mathbb R^3$ with $S_1=\{1,2\}$, $S_2=\{(1,2)\}$, $[f_1,f_2](x)=\big(\sin (x_3),-\cos (x_3), 0\big)^\top$. 
Thus, the conditions of Theorem~\ref{main} hold with $r=+\infty$, $\mu =1$. For  stabilizing system~\eqref{ex_uni} to
 a given curve $\gamma(t)\in\mathbb R^3$, we take  controls~\eqref{cont} with $k_{12}=1$:
 \begin{align}
   & u_1(t,x,\gamma)=a_1(x,\gamma)+\sqrt{\frac{4\pi|a_{12}(x,\gamma)|}{\varepsilon}}\cos\frac{2\pi t}{\varepsilon},\label{cont_uni}\\
   & u_2(t,x,\gamma)=a_2(x,\gamma)+{\rm sign}(a_{12}(x))\sqrt{\frac{4\pi|a_{12}(x,\gamma)|}{\varepsilon}}\sin\frac{2\pi t}{\varepsilon},\nonumber
 \end{align}
$$
\left(
  \begin{array}{c}
    a_1(x,\gamma) \\
    a_2(x,\gamma) \\
    a_{12}(x,\gamma) \\
  \end{array}
\right)
=\left(\begin{array}{c}
                                         (x_1-\gamma_1)\cos x_3+(x_2-\gamma_2)\sin x_3 \\
                                          x_3-\gamma_3 \\
                                          (x_1-\gamma_1)\sin x_3-(x_2-\gamma_2)\cos x_3 \\
                                        \end{array}
                                      \right).
$$
Fig.~1 (left) shows the trajectory plots of system~\eqref{ex_uni} with the curve
$
\gamma^{(1)}(t)=\big(2\cos\tfrac{t}{2}\cos t,2\cos\tfrac{t}{2}\sin t,\cos\tfrac{t}{10} \big)^\top.
$

To illustrate Corollary~1, consider the  curve
$
\gamma^{(2)}(t)=\big(3-e^{1-t},\,e^{-t^2},0 \big)^\top,
$
for which $\|\dot \gamma^{(2)}(t)\|\to 0 $ as $\to\infty$. Consequently,  $\|x(t)-\gamma^{(2)}(t)\|\to 0$ as $t\to\infty$, see Fig.~1 (middle).
\begin{remark}
The above $\gamma^{(1)}$ and $\gamma^{(2)}$ are non-admissible for system~\eqref{ex_uni}, which yields an oscillatory behavior.
Note that the asymptotic stability can be achieved for admissible curves. To illustrate this, consider the trajectory $\gamma^{(3)}(t)$ governed by
$
\dot\gamma_1^{(3)}=\dot\gamma_1^{(1)}$, $\dot\gamma_2^{(3)}=\dot\gamma_2^{(1)}$, $\dot\gamma_3^{(3)}=\frac{\dot\gamma_1^{(1)}\ddot\gamma_2^{(1)}-\dot\gamma_2^{(1)}\ddot\gamma_1^{(1)}}{{\gamma_1^{(1)}}^2+{\gamma_2^{(1)}}^2}.
$ The corresponding plot is shown in Fig.~1 (right).
\end{remark}

\subsection{Underwater vehicle}
The next example is given by the equations of motion of an autonomous 3D
underwater vehicle (see, e.g.,~\cite{Bara}):
\begin{equation}\label{ex_under}
 \dot x=\sum_{i=1}^4f_i(x)u_i,\quad x\in \mathbb R^6,\;u \in \mathbb R^4,
 \end{equation}
 where   $(x_1, x_2, x_3)$  are the coordinates of the center of mass, $(x_4$, $x_5$, $x_6)$ describe the vehicle orientation (Euler angles),  $u_1$ is the translational velocity  along the $Ox_1$ axis, and $(u_2,u_3,u_4)$ are the angular velocity components,
$$
 \begin{aligned}
&f_1(x)=(\cos x_5\cos x_6, \cos x_5\sin x_6,{-}\sin x_5,0,0,0)^\top, \\
&f_2(x) =(0,0,0,1,0,0)^\top,
\end{aligned}
$$
$$
 \begin{aligned}
 & f_3(x){=}(0,0,0,\sin x_4{\rm tg}\, x_5,\cos x_4,\sin x_4\sec x_5)^\top,\\
 &f_4(x){=}(0,0,0,\cos x_4{\rm tg}\, x_5,{-}\sin x_4,\cos x_4\sec x_5)^\top.
 \end{aligned}
$$
The rank condition~\eqref{rank} is satisfied in $D =\{x \in R^6:-\tfrac{\pi}{2}<x_5<\tfrac{\pi}{2}\}$ with $S=\{(1,3),(1,4)\}$. Therefore, the matrix
$$
\mathcal F(x)=\left(
f_1(x),\ f_2(x),\ f_3(x),\ f_4(x),\ [f_1,f_3](x),\ [f_1,f_4](x)
              \right)
$$
is nonsingular in $D$.
Thus, controls~\eqref{cont} take the form
 \begin{align}
    u_1(t,x,\gamma)=&a_1(x,\gamma)+\sqrt{\frac{4\pi|a_{13}(x,\gamma)|}{\varepsilon}}\cos\frac{2\pi k_{13}t}{\varepsilon}\nonumber\\
   &+\sqrt{\frac{4\pi|a_{14}(x,\gamma)|}{\varepsilon}}\cos\frac{2\pi k_{14}t}{\varepsilon},\nonumber\\
   u_2(t,x,\gamma)=&a_2(x,\gamma),\label{cont_under}\\
    u_3(t,x,\gamma)=&a_3(x,\gamma)+{\rm sign}(a_{13}(x))\sqrt{\frac{4\pi|a_{13}(x,\gamma)|}{\varepsilon}}\sin\frac{2\pi k_{13} t}{\varepsilon},\nonumber\\
    u_4(t,x,\gamma)=&a_4(x,\gamma)+{\rm sign}(a_{14}(x))\sqrt{\frac{4\pi|a_{14}(x,\gamma)|}{\varepsilon}}\sin\frac{2\pi k_{14} t}{\varepsilon},\nonumber
 \end{align}
with
$
a(x,\gamma)=-\alpha \mathcal F^{-1}(x)(x-\gamma).
$

For the illustration, take
$
\gamma^{(4)}(t)=\left(\cos\tfrac{t}{4},\,\tfrac{t}{4},\,\sin\tfrac{t}{4},\,0,\,0\,0\right)^\top.
$
The results of numerical simulations are shown in Fig.~2. Note that the curve $\gamma^{(4)}(t)$ is non-admissible for system~\eqref{ex_under}, which results in an oscillatory behavior of the trajectories.
\begin{figure*}[h!]
 \begin{minipage}{0.5\linewidth}
  \includegraphics[width=1\linewidth,height=0.6\linewidth]{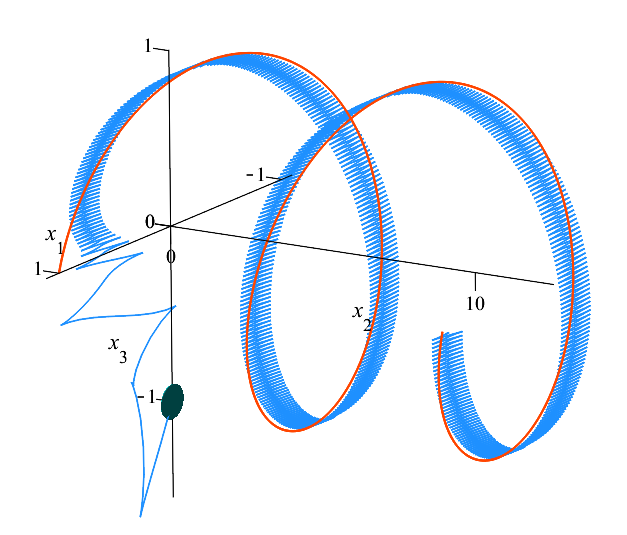}
 \end{minipage}\hfill
  \begin{minipage}{0.5\linewidth}
\includegraphics[width=1\linewidth]{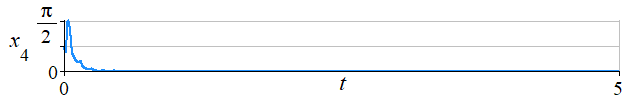}
\includegraphics[width=1\linewidth]{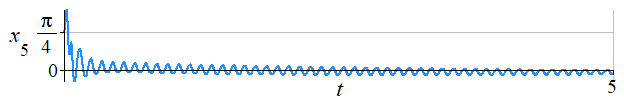}
\includegraphics[width=1\linewidth]{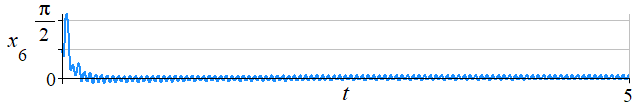}
\end{minipage}
\caption{Trajectories of system~\eqref{ex_under} with controls~\eqref{cont_under}; $\alpha=15$, $\varepsilon=0.1$, $k_{13}=1$, $k_{14}=2$, $x^0=(0,0,-1,\tfrac{\pi}{4},\tfrac{\pi}{4},\tfrac{\pi}{4})^\top$.}
\end{figure*}
\subsection{Rear-wheel driving car}
The proposed approach can also be extended to nonholonomic systems of higher degrees. For systems of degree two, it is possible to use a control design scheme similar to that introduced in~\cite{GZ18,ZG17}.
For example, consider a kinematic model of a rear-wheel driving car proposed in~\cite{Lu98}:
\begin{equation}\label{ex_car}
  \dot x= f_1(x)u_1+f_2(x)u_2,\quad x\in\mathbb R^4,\,u\in\mathbb R^2,
\end{equation}
where $(x_1,x_2)$ are the
Cartesian coordinates of the rear wheel, $x_3$ is the steering angle, $x_4$ specifies the orientation of the car
body with respect to the $x_1$ axis, $u_1$ and $u_2$ are the driving and the steering velocity input, respectively,
$$
\begin{aligned}
&f_1(x)=(\cos x_4,\,\sin x_4,\,0,\tan x_3)^\top, \ f_2(x) =(0,0,1,0)^\top.\\
 \end{aligned}
$$
In this case,
$
{\rm span}\{f_1(x),\,f_2(x),\,[f_1,f_2](x),\,[[f_1,f_2],f_1](x)\}=\mathbb R^4
$
for all $x\in D =\{x \in R^4:-\tfrac{\pi}{2}<x_3<\tfrac{\pi}{2}\}$. Following the control design scheme from~\cite{GZ18}, we take
 \begin{align}
 u_1&(t,x,\gamma)=a_1(x,\gamma)+\sqrt{\frac{4\pi|a_{12}(x,\gamma)|}{\varepsilon}}\cos\frac{2\pi k_{12}t}{\varepsilon}\nonumber\\
   &+\sqrt[3]{\frac{16\pi^2(k_2^2-k_1^2)a_{121}(x,\gamma)}{\varepsilon^2}}\cos\frac{2\pi k_{1}t}{\varepsilon}\Big(1+\sin\frac{2\pi k_{2}t}{\varepsilon}\Big),\nonumber\\
 u_2&(t,x,\gamma)=a_2(x,\gamma)+{\rm sign}(a_{12}(x))\sqrt{\frac{4\pi|a_{12}(x,\gamma)|}{\varepsilon}}\sin\frac{2\pi k_{12} t}{\varepsilon}  \nonumber \\ &+\sqrt[3]{\frac{16\pi^2(k_2^2-k_1^2)a_{121}(x,\gamma)}{\varepsilon^2}}\sin\frac{2\pi k_{2}t}{\varepsilon},\label{cont_car}
 \end{align}
with the vector of coefficients
$
a(x,\gamma)=-\alpha \mathcal F^{-1}(x)(x-\gamma)
$ and
$
\mathcal F(x)=\left(
f_1(x)\ f_2(x)\  [f_1,f_2](x),\ \big[[f_1,f_2],f_1\big](x)
              \right).
$

Fig.~3 presents the trajectory plots of system~\eqref{ex_car}--\eqref{cont_car} for a non-admissible curve
$
\gamma^{(4)}(t)=\left(5\sin\tfrac{t}{4}, 5\sin\tfrac{t}{4} \cos\tfrac{t}{4},\,0,\,0\right)^\top.
$
\begin{figure}[hb]
 \begin{minipage}{1\linewidth}
  \includegraphics[width=1\linewidth]{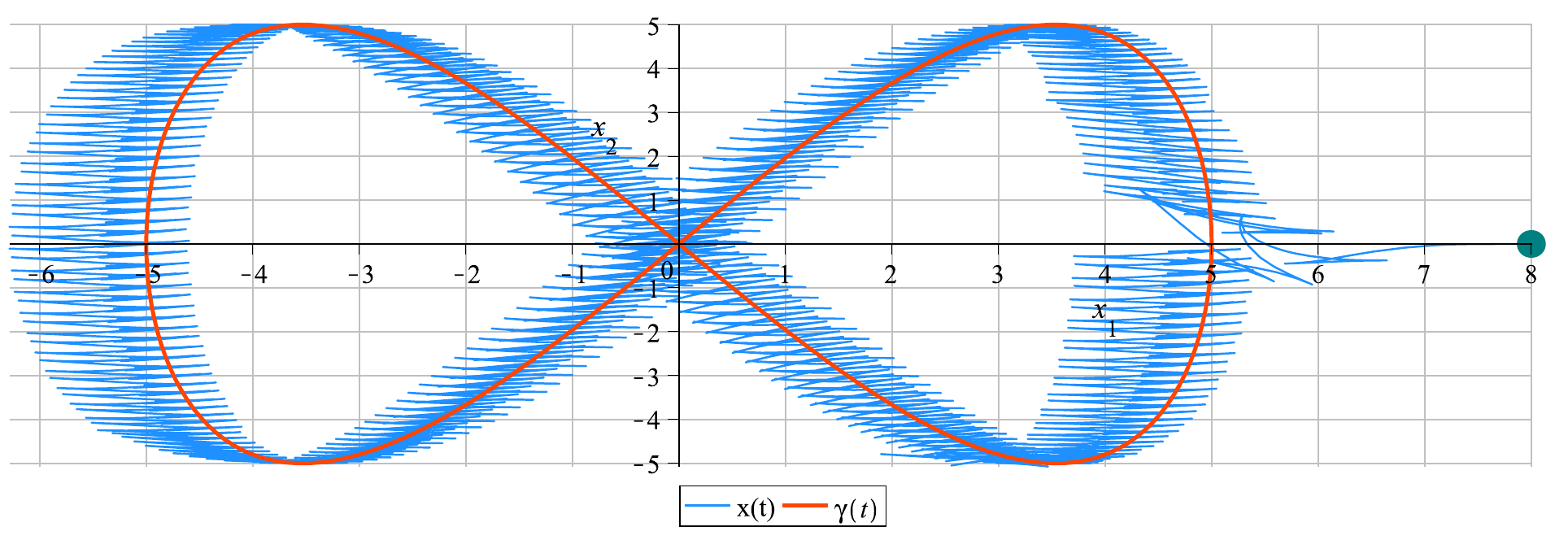}
\includegraphics[width=1\linewidth]{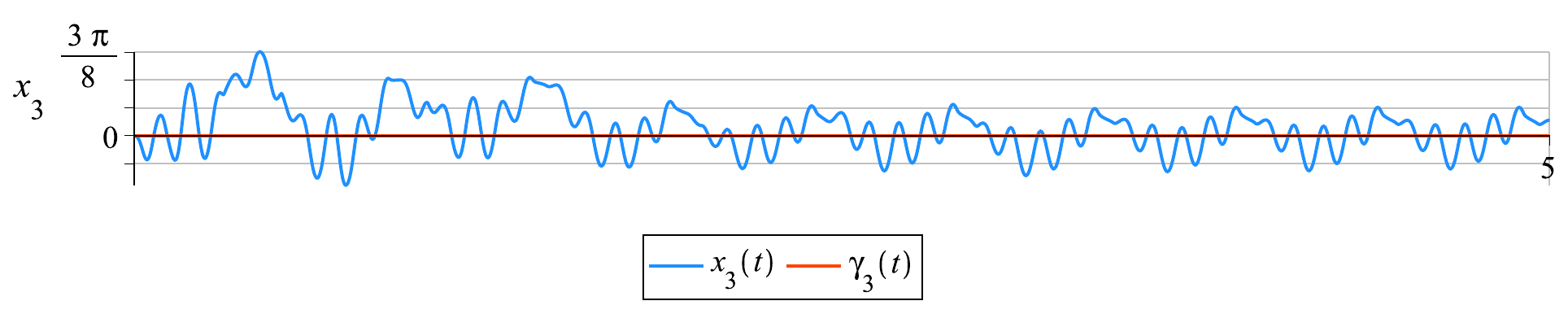}
\includegraphics[width=1\linewidth]{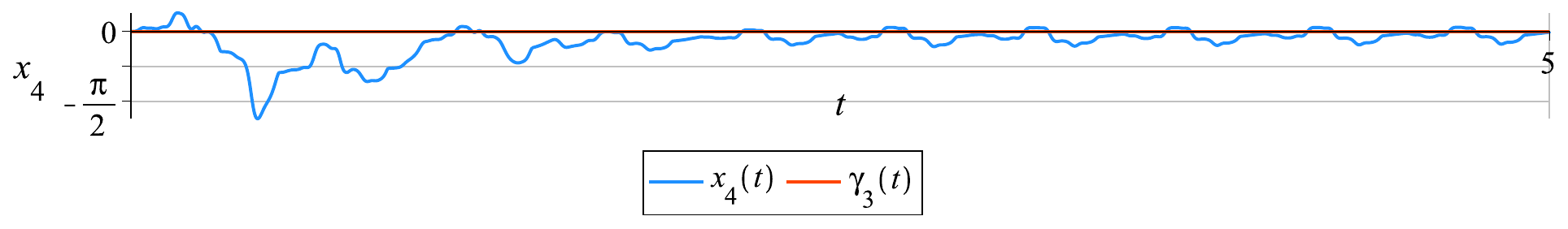}
\caption{Trajectories of system~\eqref{ex_car} with controls~\eqref{cont_car}; $\alpha=5$, $\varepsilon=0.5$, $x^0=(8,0,0,0)^\top$.}
 \end{minipage}\hfill
\end{figure}

\section{CONCLUSIONS AND FUTURE WORK}
The above numerical simulations confirm that the proposed controller~\eqref{cont} can be used for approximate tracking of reference curves under an appropriate choice of parameters $\alpha$ and $\varepsilon$.
By comparing the left and right plots in Fig.~1, we note that the amplitude of oscillations near non-admissible curve (Fig.~1, left) significantly exceeds the deviation from the admissible curve (Fig.~1, right).
This feature underlines the assertion of Corollary~1 and illustrates the essence of our approach for considering the stability of a family of sets.
The example in Section~III.C shows that our approach can also be extended to nonholonomic systems of higher degrees.
We do not study here the stabilization problem under general controllability conditions, leaving this issue for future work.

\appendix
\section{Proof of the main result}
\subsection{Proof of Theorem~\ref{main}}
To prove Theorem~1, we
 will use the following result.
\begin{lemma}[\cite{Zu16}]~\label{lemma_x}
 \emph{Let $\tilde D\subseteq\mathbb R^n$ be a convex domain, and let $x(t)\in \tilde D$, $0\le t\le\tau$, be a  solution of system~\eqref{Sigma} with some control $u\in C[0,\tau]$. Assume that there exist  $M ,L>0$ such that
$
\|f_i(x)\|\le   M,\,\|f_i(x)-f_i(y)\|\le  L\|x-y\|,
$
for all $x,y\in \tilde D$, $i{=}\overline{1,m}$.
    Then
    \begin{equation}\label{est_Z}
      \|x(t)-x(0)\|{\le}\tfrac{ M}{ L}(e^{U Lt}{-}1),\;t{\in}[0,\tau],
    \end{equation}
    with $\displaystyle U=\max\limits_{t\in[0,\tau]}\sum_{i=1}^{m}|u_{i}(t)|$.}
  \end{lemma}
  \begin{lemma}[\cite{Lam96,GZE18}]\label{volterra}
\emph{Let the vector fields $f_i$ be Lipschitz continuous in a domain $D\subseteq\mathbb R^n$, and $f_i\in C^2(D\setminus\Xi)$, where $\Xi=\{x\in D:f_i(x)=0\text{ for all }1\le i\le m\}$. Assume, moreover, that ${ L}_{f_j}f_i, { L}_{f_l}{ L}_{f_j}f_i\in
C(D;\mathbb R^n)$,   for all $i,j,l=\overline{1,m}$. If $x(t)\in D$, $t\in[0,\tau]$, is a  solution of system~\eqref{Sigma} with $u\in C[0,\tau]$ and $x(0)=x^0\in D$, then  $x(t)$ can be represented by the Volterra series:
\begin{align}
   & x(t){=}x^0{+}{\sum_{i=1}^{m}}f_i(x^0)\int\limits_0^t u_{i}(v)dv \label{volt1}\\
   &{+}\sum_{\hspace{-0.75em}i,j=1}^{m} L_{f_j}f_i(x^0)\int\limits_0^t\int\limits_0^v  u_{i}(v) u_{j}(s)dsdv+R(t),\,t\in[0,\tau], \nonumber
\end{align}
where
$$
R(t){=}{\sum\limits_{\hspace{-0.75em}i,j,l=1}^{m}}{\int\limits_0^t}{\int\limits_0^v}{\int\limits_0^s}{ L}_{f_l}{ L}_{f_j}f_i(x(p))  u_{i}(v)u_{j}(s)u_{l}(p)\,dp\,ds\,dv \nonumber
$$
     is the remainder of the Volterra series expansion.}
     \end{lemma}
  \subsection*{Proof of Theorem~\ref{main}}

  Let us take any positive numbers $\delta$, $\delta'$, and $\rho'$ from the inequalities
  $
\frac{\nu}{\alpha}<\rho'< \rho < \delta < \delta' < r,
$
and denote $\displaystyle D'=\bigcup_{t\ge 0}B_{\delta'}(\gamma(t))$,  $\gamma^0=\gamma(0)$.
It is clear that
 $$B_\delta(\gamma(\tau))\subset  D' \subset D\quad \text{for each}\; \tau\ge 0.
  $$

Let $x^0\in B_\delta(\gamma^0)$.
  Our first goal is to find an $\varepsilon_1>0$ such that the corresponding solution of system~\eqref{Sigma} with the initial condition $x^0\in \overline{B_\delta(\gamma^0)}$ and controls~\eqref{cont}
  is well-defined on $[0,\varepsilon]$ and satisfies the property $x(t)\in D$ for all $t\in [0,\varepsilon]$, $\varepsilon\in(0,\varepsilon_1)$.
 Let the control functions be defined by~\eqref{cont}, and let
$
M_1=\sup\limits_{\hspace{-1em}\underset{1\le i\le m}{x\in D'}}\|f_i(x)\|,\;U(x^0,\gamma^0)=\max\limits_{0\le t\le \varepsilon}\sum_{i=1}^m |u_i^\varepsilon(t,x^0,\gamma^0)|.
$
Using H\"{o}lder's inequality, one can estimate $U(x^0,\gamma^0)$ as
\begin{align}
U(x^0,\gamma^0)
\le C_1\|x^0-\gamma^0\|+\frac{C_2}{\sqrt\varepsilon}\sqrt{\|x^0-\gamma^0\|},\label{est_u}
\end{align}
with $C_1= \alpha\mu \sqrt{|S_1|}$, $C_2=4\sqrt{{\pi\mu\alpha}}\Big(\sum_{(j_1,j_2)\in S_2}{\kappa_{j_1j_2}}^{2/3}\Big)^{3/4}.$
Note also that
\begin{equation}\label{geps}
\|\gamma(t)-\gamma^0\|= \left\|\int_0^t \dot\gamma(s)ds\right\|\le \nu t \quad \text{for all}\; t\in [0,\varepsilon].
\end{equation}
From Lemma~\ref{lemma_x} and estimate~\eqref{est_u},
$$
\begin{aligned}
\|x(t)-x^0\|&\le \frac{M_1}{L}(e^{U(x^0,\gamma^0) L\varepsilon}{-}1)\\
&\le  \frac{M_1}{L}(e^{L( C_1\varepsilon\delta'+{C_2}\sqrt{\varepsilon\delta'})}{-}1),\text{ for all }t{\in}[0,\varepsilon].
\end{aligned}
$$
Thus, defining $d=\min\{\delta'-\delta,\tfrac{\rho-\rho'}{2}\}>0$ and
$$
\begin{aligned}
\varepsilon_{1}=\min\Big\{&\frac{\rho-\rho'}{2\nu},\\
&\frac{1}{4\delta'}\Big(\sqrt{\Big(\frac{C_2}{C_1}\Big)^2+\frac{4}{L C_1}\ln\Big(\frac{dL}{M_1}+1\Big)}-\frac{C_2}{C_1}\Big)^2\Big\},
\end{aligned}
$$
we conclude that   for any $\varepsilon\in(0,\varepsilon_1)$,
\begin{equation}\label{xtx0}
\|x(t)-x^0\|<d\text{ for all }t\in[0,\varepsilon],
\end{equation}
so that the solutions of system~\eqref{Sigma} with controls~\eqref{cont}  and the initial conditions $x(0)=x^0\in \overline{B_\delta(0)}$ stay in $D$ for all $t\in[0,\varepsilon]$.
Moreover,~\eqref{geps} and~\eqref{xtx0} yield that  if $x^0\in B_{\rho'}(\gamma^0)$, then $x(t)\in \mathcal S_t^\rho$ for each $t\in[0,\varepsilon]$:
 \begin{equation}\label{xtgt}
\begin{aligned}
 \|x(t)-\gamma(t)\|&\le\|x(t)-x^0\|+\|\gamma(t)-\gamma^0\|+\|x^0-\gamma^0\|\\
 &<d+\varepsilon\nu+\rho'<\rho.
\end{aligned}
 \end{equation}
 Using Lemma~\ref{volterra} we obtain the following representation of the solutions of system~\eqref{Sigma} with  the controls defined by~\eqref{cont} and the initial conditions $x(0)=x^0\in D$:
 \begin{align}\label{xeps}
x(\varepsilon)=x^0-  \varepsilon \alpha(x^0-\gamma^0)+ R(\tilde a,\varepsilon),
\end{align}
where
\begin{align*}
 R( a,&\varepsilon)={\varepsilon^{3/2}}\sum_{j_1\in S_1}\sum_{j_2=1}^m[f_{j_1},f_{j_2}](x^0)a_{j_1}(x^0,\gamma^0)\\
 &\qquad\times\sum_{q:(q,j_2)\in S_2}{\rm sign}(a_{qj_2}(x^0,\gamma^0))\sqrt{\frac{|a_{qj_2}(x^0,\gamma^0)|}{\pi \kappa_{qj_2}}}\\
&+\frac{\varepsilon^2}{2}\sum_{j_1,j_2\in S_1} L_{f_{j_2}}f_{j_1}(x^0)a_{j_1}(x^0,\gamma^0)a_{j_2}(x^0,\gamma^0)+\tilde r(\varepsilon).
  \end{align*}
  Denote
  $$M_2{=}\sup\limits_{\hspace{-1em}\underset{j_1,j_2=\overline{1,m}}{x\in D'}}\| L_{f_{j_1}}f_{j_2}(x)\|,
  M_3{=}\frac{1}{6}\max\limits_{x\in D'}{\sum\limits_{j_1,j_2,j_3=1}^m}\big\| L_{f_{j_3}} L_{f_{j_2}}f_{j_1}(x)\big\|.
  $$
   Then from~\eqref{est_u},
\begin{equation*}\label{r_est}
   \begin{aligned}
     \|\tilde r(\varepsilon)\|&\le {M_3}\big(U(x^0,\gamma^0)\varepsilon\big)^3\\
     &\le {M_3}\Big(C_1\varepsilon \|x^0-\gamma^0\|
     +{C_2}\sqrt{\varepsilon\|x^0-\gamma^0\|}\Big)^3
   \end{aligned}
\end{equation*}
for all $t\in[0,\varepsilon]$,  and
  \begin{align*}
  \| R( a,\varepsilon)\|&\le 2\varepsilon^{3/2}\|a(x^0,\gamma^0)\|^{3/2}M_2\sqrt{|S_1|}\\
  &\times\sum_{j_1=1}^m\Big(\sum_{(j_2,j_1)\in S_2}\kappa_{j_2j_1}^{-2/3}\Big)^{3/4}+\varepsilon^2\|a(x^0,\gamma^0)\|^2\frac{M_2}{2}\\
  &+  M_3\Big(C_1\varepsilon \|x^0-\gamma^0\|+{C_2}{\sqrt\varepsilon}\sqrt{\|x^0-\gamma^0\|}\Big)^3\\
  &\le \sigma \varepsilon^{3/2}\|x^0-\gamma^0\|^{3/2},
\end{align*}
where
\begin{align*}
\sigma =M_2\bigg(&2\left(\alpha\mu\right)^{3/2}\sqrt{|S_1|}\sum_{j_1=1}^m\Big(\sum_{(j_2,j_1)\in S_2}\kappa_{j_2j_1}^{-2/3}\Big)^{3/4}\\
&+\frac{1}{2}\sqrt{\varepsilon \delta'}(\alpha\mu)^2\bigg)+M_3\left(C_2+C_1\sqrt{\varepsilon \delta'}\right)^3.
\end{align*}

Considering~\eqref{xeps} and~\eqref{geps}, we obtain
$$
\|x(\varepsilon)-\gamma(\varepsilon)\|\le\|x^0-\gamma^0\|\big(1-\varepsilon(\alpha-\sigma \sqrt{\varepsilon\|x^0-\gamma^0\|})\big)+\varepsilon\nu.
$$
Recall that $\rho'\in\Big(\frac{\nu}{\alpha},\rho\Big)$.
For any $\lambda\in(0,\alpha-\tfrac{\nu}{\rho'})$, let $\varepsilon_2=\frac{1}{\delta'}\Big(\frac{\alpha-\lambda}{\sigma}-\frac{\nu}{\sigma\rho'}\Big)^2,\Big(\lambda+\frac{\nu}{\rho'}\Big)^{-1}.$
 Then, for any $\varepsilon\in(0,\min\{\varepsilon_1,\varepsilon_2\})$, $x^0\in D'$,
\begin{equation}\label{xege}
\|x(\varepsilon)-\gamma(\varepsilon)\|\le\|x^0-\gamma^0\|\Big(1-\varepsilon\Big(\lambda+\frac{\nu}{\rho'}\Big)\Big)+\varepsilon\nu.
\end{equation}
Consider two cases.

Case 1) If $x^0\in \mathcal S_{0}^{\rho'}$, then it is easy to see from~\eqref{xege} that $x(\varepsilon)\in \mathcal S_{0}^{\rho'}$. From~\eqref{xtgt}, $x(t)\in \mathcal S_t^\rho$ for each $t\in[0,\varepsilon]$.

Case 2) Assume now that $\|x^0-\gamma^0\|>\rho'$. Then
$$
\begin{aligned}
\|x(\varepsilon)-\gamma(\varepsilon)\|&\le\|x^0-\gamma^0\|\Big(1-\varepsilon\Big(\lambda+\frac{\nu}{\rho'}-\frac{\nu}{\|x^0-\gamma^0\|}\Big)\Big) \\
&<\|x^0-\gamma^0\|\Big(1-\varepsilon\lambda\Big).
\end{aligned}
$$
Iterating the above inequality for $x(t_0)\in\mathcal S_{t_0}^\delta$,  we conclude that there exists an $N\in\mathbb N$ such that $\|x(t)-\gamma(t)\|>\rho$ for each $t=0,\varepsilon,2\varepsilon,\dots,(N-1)\varepsilon$, and $x(N\varepsilon)\in\mathcal S_{N\varepsilon}^{\rho'}$ (this can be proved by contradiction).
Repeating the argumentation of Case 1) and Case 2), we conclude that $x(t)\in\mathcal S_{t}^{\rho}$ for all $t\ge N\varepsilon$.

It remains to consider an arbitrary $t\in[0,N\varepsilon]$. Denote by $t_{in}=\Big[\frac{t}{\varepsilon}\Big]$ the integer part of $\frac{t}{\varepsilon}$. Since $t-t_{in}\varepsilon<\varepsilon$, we have
$$
\begin{aligned}
\|x(t)-\gamma(t)\|&\le \|x(t_{in}\varepsilon)-\gamma(t_{in}\varepsilon)\| + \|x(t)-x(t_{in}\varepsilon)\| \\
&+\|\gamma(t)-\gamma(t_{in}\varepsilon)\| \\
&\le\|x^0-\gamma^0\|e^{-\lambda t_{in}\varepsilon}+\frac{M_1}{L}(e^{LU(t_{in}\varepsilon)\varepsilon}{-}1)+\varepsilon\nu,
\end{aligned}
$$
where $U(t_{in}\varepsilon)=\max\limits_{s\in[ t_{in}\varepsilon,t]}\sum_{i=1}^{m}|u^{\varepsilon}_{i}(s,x(t_{in}\varepsilon),\gamma(t_{in}\varepsilon))|.$
From~\eqref{est_u},
\begin{align*}
U(t_{in}\varepsilon)&\le C_1\|x(t_{in}\varepsilon)-\gamma(t_{in}\varepsilon)\|+\frac{C_2}{\sqrt\varepsilon}\sqrt{\|x(t_{in}\varepsilon)-\gamma(t_{in}\varepsilon)\|}\\
&\le C_1\|x^0-\gamma^0\|e^{-\lambda t_{in}\varepsilon}+\frac{C_2}{\sqrt\varepsilon}\sqrt{\|x^0-\gamma^0\|e^{-\lambda t_{in}\varepsilon}}.
\end{align*}
Let $\varepsilon_3\in\Big(0,\frac{1}{\delta'}\Big(\sqrt{\frac{C_2^2}{4C_1^2}+\frac{1}{L}}-\frac{C_2}{2C_1}\Big)^2\Big)$. Then, for any $\varepsilon\in(0,\min\{\varepsilon_1,\varepsilon_2,\varepsilon_3\})$,
$$
\begin{aligned}
\|&x(t)-\gamma(t)\|\le \|x^0-\gamma^0\|e^{-\lambda t_{in}\varepsilon}+{M_1}(e-1)U(t_{in}\varepsilon)+\varepsilon\nu\\
& \le\|x^0-\gamma^0\|e^{-\lambda t_{in}\varepsilon}+M_1(e-1)\Big(C_1\varepsilon\|x^0-\gamma^0\|e^{-\lambda t_{in}\varepsilon}\\
&+{C_2}{\sqrt\varepsilon}\sqrt{\|x^0-\gamma^0\|e^{-\lambda t_{in}\varepsilon}}\Big)+\varepsilon\nu\\
&\le \kappa(\|x^0-\gamma^0\|) e^{-\lambda_1 t}+\varepsilon\nu<\kappa e^{-\lambda t_{in}\varepsilon}+\frac{\rho-\rho'}{2},
\end{aligned}
$$
where $\lambda_1=\lambda/2$,
$$
\kappa(\|x^0-\gamma^0\|)=(1+\varepsilon M_1(e-1)C_1)\|x^0-\gamma^0\|e^{\varepsilon}+{C_2}{\sqrt\varepsilon}\sqrt{\|x^0-\gamma^0\|}e^{\varepsilon/2}.
$$

{Thus, for any $x^0\in\mathbb \mathcal S^{\delta}_0$, there exists a $T\ge 0$ such that  ${\rm dist}(x(t), \mathcal S_t^{\rho'}) \le \kappa(\|x^0-\gamma^0\|) e^{-\lambda_1 t}$ for all $t\in[0,T)$, and $x(t)\in\mathcal S_t^{\rho}$ for all $t\ge T$, which proves Theorem~1.}


\begin{thebibliography}{99}
\bibitem{Ai05} A. Ailon, N. Berman, and S. Arogeti, On controllability and trajectory tracking of a kinematic vehicle model, Automatica, Vol.~41, no.~5, pp.~889--896, 2005.
\bibitem{Ali16} Z.A. Ali, D. Wang, M. Safwan, W. Jiang, and M. Shafiq, Trajectory Tracking of a Nonholonomic Wheeleed Mobile Robot Using Hybrid Controller, International Journal of Modeling and Optimization, Vol.~6, No.~3, pp.~136--141, 2016.
\bibitem{Dandr95} B. d'Andr\'{e}a-Novel, G. Campion, and G. Bastin, Control of nonholonomic wheeled mobile robots by state feedback linearization, The International Journal of Robotics Research, Vol.~14, No.~6, pp.~543--559, 1995.
\bibitem{Ast94} A. Astolfi, On the stabilization of nonholonomic systems, in Proc. 33rd IEEE Conf. on Decision and Control, Vol.~4, pp.~3481--3486, 1994.
\bibitem{Bara} {J.~Barraquand and J.-C.~Latombe}, On non-holonomic mobile robots and optimal maneuvering,Revue d'Intelligence Artificielle, Vol.~13, pp.~77--103, 1989.
\bibitem{Bloch92} A.\,M. Bloch,  M. Reyhanoglu, and N.\,H. McClamroch, Control and stabilization of nonholonomic dynamic systems, IEEE Transactions on Automatic control, Vol.~37, No.~11, pp.~1746--1757, 1992.
\bibitem{Bro81} R.\,W. Brockett, Asymptotic stability and feedback stabilization, in Differential Geometry Control Theory, R.W. Brockett, R.S. Hillman, H.J. Sussmann, Eds.  Boston: Birkh{\"a}user, 1983, pp.~181--191.
\bibitem{Cor92} J.-M. Coron, Global asymptotic stabilization for
controllable systems without drift, Math. Control Signals Systems, vol.~5, pp.~295--312, 1992.
\bibitem{Lu98} A. De Luca, G. Oriolo, and C. Samson, Feedback control of a nonholonomic car-like robot, In Robot Motion Planning and Control, J.-P. Laumond, Ed. Berlin, Heidelberg: Springer,  1998, pp.~171--253.
\bibitem{Dong99} W. Dong, W.\,L. Xu, and W. Huo, Trajectory tracking control of dynamic non‐holonomic systems with unknown dynamics, International Journal of Robust and Nonlinear Control, Vol.~9, pp.~905--922, 1999.
\bibitem{Fl95} M. Fliess, J. L\'{e}vine, P. Martin, and P. Rouchon, Design of trajectory stabilizing feedback for driftless at systems, in Proc. European Control Conf. 1995, pp.~1882--1887.
\bibitem{GDEZ17} V. Grushkovskaya, H.-B. D\"{u}rr, C. Ebenbauer, and A. Zuyev, Extremum seeking for time-varying functions using Lie bracket approximations, IFAC-PapersOnLine, vol.~50, pp.~5522--5528, 2017.
\bibitem{GMZME18} V. Grushkovskaya, S. Michalowsky, A. Zuyev, M. May, and C. Ebenbauer, A family of extremum seeking laws for a unicycle model with a moving target: theoretical and experimental studies, in Proc. 18th European Control Conf, pp.~912--917, 2018.
\bibitem{GZ18} V. Grushkovskaya and  A. Zuyev, Obstacle avoidance problem for second degree nonholonomic systems, in Proc. 57th IEEE Conf. on Desicion and Control, pp.~1500--1505, 2018.
\bibitem{GZE18} V. Grushkovskaya, A. Zuyev, and C.~Ebenbauer, On a class of generating vector fields for the extremum seeking problem: Lie bracket approximation and stability properties, Automatica, Vol.~94, pp.~151--160, 2018.
\bibitem{Khalil}H.K.~Khalil, Nonlinear Systems. 3rd Ed., Prentice Hall, 2002.
\bibitem{Kol95} I. Kolmanovsky and N.\,H. McClamroch, Developments in nonholonomic control problems, IEEE Control Systems, Vol.~15, No.~6, pp.~20--36, 1995.
 \bibitem{Lam96} F.Lamnabhi-Lagarrigue, {V}olterra and {F}liess series expansions for nonlinear systems, in The Control Handbook, W.~S. Levine, ed. CRC press, 1996,  pp.~879--888.
\bibitem{La02} J.\,A.  Langa, J.\,C. Robinson  and  A. Su\'{a}rez, Stability, instability, and bifurcation phenomena in non-autonomous differential equations, Nonlinearity, Vol.~15, No.~3, pp.~887--903, 2002.
\bibitem{Mag17} M. Maghenem, A. Loria, and E. Panteley, Global tracking-stabilization control of mobile robots with parametric uncertainty,  IFAC-PapersOnLine, Vol.~50, No.~1, pp.~4114--4119, 2017.
\bibitem{Pom92} P. Morin, J.-B. Pomet, and C. Samson,
Design of homogeneous time-varying stabilizing control laws for driftless controllable systems via oscillatory approximation of Lie brackets in closed loop,
SIAM J. Control Optim., vol.~38, pp.~22--49, 1999.
\bibitem{MS08b} P. Morin and C. Samson, Trajectory tracking for nonholonomic systems. Theoretical background
and applications. [Research Report] 2008, 49 p. $<$inria-00260694v1$>$.
\bibitem{Pan11} D. Panagou, H.\,G. Tanner, and K.\,J. Kyriakopoulos,  Control of nonholonomic systems using reference vector fields. In Proc. 50th IEEE Conf. on Decision and Control and European Control Conference, pp.~2831--2836, 2011.
%
\bibitem{Rav18} H. Ravanbakhsh, S. Aghli, C. Heckman, and S. Sankaranarayanan, Path-Following through Control Funnel Functions, arXiv preprint, arXiv:1804.05288, 9 p., 2018.
\bibitem{Sar17} M. Sarfraz and  F. Rehman, Feedback Stabilization of Nonholonomic Drift-Free Systems Using
Adaptive Integral Sliding Mode Control, Arabian Journal for Science and Engineering, Vol.~42, pp.~2787--2797, 2017.
\bibitem{Tian02} Y.-P. Tian and S. Li, Exponential stabilization of nonholonomic dynamic systems by smooth time-varying control, Automatica, Vol. 38, No. 7, pp.~1139--1146, 2002.
\bibitem{Walsh92} G. Walsh, D. Tilbury, S. Sastry, R. Murray, and J.\,P. Laumond, Stabilization of trajectories for systems with nonholonomic constraints, in Proc. IEEE International Conference on Robotics and Automation, pp.~1999--2004, 1992.
\bibitem{Wang15} Y. Wang, Z. Miao, H. Zhong, and Q. Pan, Simultaneous stabilization and tracking of nonholonomic mobile robots: A Lyapunov-based approach, IEEE Transactions on Control Systems Technology, Vol. 23, No. 4, pp.~1440--1450, 2015.
\bibitem{Yu15} X. Yu, L. Liu, and G. Feng, Trajectory Tracking for Nonholonomic Vehicles with Velocity Constraints, IFAC-PapersOnLine, vol. 48, no. 11, pp.~918--923, 2015.
\bibitem{Zu16} A. Zuyev, Exponential stabilization of nonholonomic systems by means of oscillating controls,  SIAM Journal on Control and Optimization, Vol.~54, no.~3, pp.~1678--1696, 2016.
\bibitem{ZG16} A. Zuyev, V. Grushkovskaya, and P.~Benner, Time-varying stabilization of a class of driftless systems satisfying second-order controllability conditions, in Proc. European Control Conf. 2016,  pp.~575--580, 2016.
\bibitem{ZG17} A. Zuyev and V. Grushkovskaya, Motion planning for control-affine systems satisfying low-order controllability conditions, International Journal of Control, vol. 90, no. 11, pp.~2517--2537, 2017.
\end{thebibliography}
\end{document}